\input{aipcheck}
\documentclass[,final ] {aipproc}
\layoutstyle{8x11single}
\usepackage{amsmath}
\usepackage{amsfonts}
\usepackage{amssymb}

\usepackage{epsfig} %% for loading postscript figures
\usepackage{graphicx}
\usepackage{enumerate}
\usepackage{{subfigure}}
\usepackage{epstopdf}
\usepackage{enumerate}
\begin{document}
\begin{center}\small{In the name of Allah, the Beneficent, the Merciful.}\end{center}
\vspace{2cm}
	
\title{A note on a separating system of rational invariants for finite dimensional generic algebras}
\author{U. Bekbaev}
{
address={ Deparment of Science in Engineering, KOE, IIUM, Malaysia.\\
\hfill \break
bekbaev@iium.edu.my \\}}

% ----------------------------------------------------------------

\begin{abstract} The paper deals with a construction of a separating system of rational invariants for finite dimensional generic algebras. In the process of dealing an approach  to a rough classification of finite dimensional algebras is offered by attaching them some quadratic forms.
\end{abstract}
\maketitle
% ----------------------------------------------------------------

\section{Introduction}

In \cite{B2017} we have offered an approach to classification problem of finite dimensional algebras  with respect to basis changes. It has been shown that if one has a special map with some properties then he is able to classify, to list canonical representations, all algebras who's set of structural constants , with respect to a fixed basis, do not nullify some polynomial. In this case he is also able to provide a separating system of rational invariants for those algebras. It was successfully applied in \cite{A2017} to get a complete classification of all $2$-dimensional algebras over algebraically closed fields. 

Unfortunately, so far we have no example of such a map in $3$-dimensional case. Therefore in the current paper we deal with a weaker problem, namely with a construction of  separating system of rational invariants for finite dimensional generic algebras. The theoretical existence of such system of invariants is known \cite{P2014}. By generic algebras we mean the set of all algebras who's system of structural constants does not nullify a fixed nonzero polynomial in structural variables, over the basic field $F$. In process of dealing with the problem we show a way for a rough classification of finite dimensional algebras by attaching them some quadratic forms. 

The next section contains the main results.

\section{Main results}

Further whenever $A=(a_{ij})\in Mat(p\times q, F)$, $B\in Mat(p'\times q', F)$ we use $A\otimes B$ for the matrix \[\left(\begin{array}{cccc} a_{11}B&a_{12}B&...&a_{1q}B\\
a_{21}B&a_{22}B&...&a_{2q}B\\
.&.&\ddot{.}&.\\
a_{p1}B&a_{p2}B&...&a_{pq}B\end{array}\right),\ \mbox{where}\ F-\ \mbox{is a field of characteristic not}\ 2. \]
Let us consider any $m$-dimensional algebra $\mathbf{A}$ with multiplication $\cdot$ given by a bilinear map $(\mathbf{u},\mathbf{v})\mapsto \mathbf{u}\cdot \mathbf{v}$. If $e=(e_1,e_2,...,e_m)$ is a
basis for $\mathbf{A}$ then one can represent the bilinear map by a matrix \[A_e=(A^i_{ejk})_{i,j,k=1,2,...,m}\in Mat(m\times m^2;F),\] where $e_j\cdot e_k=e_1A^1_{ejk}+e_2A^2_{ejk}+...+e_mA^m_{ejk}$, $j,k=1,2,...,m$, such that \[\mathbf{u}\cdot \mathbf{v}=eA_e(u\otimes v)\] for any $\mathbf{u}=eu,\mathbf{v}=ev,$
where $u=(u_1, u_2,..., u_m), v=(v_1, v_2,..., v_m)$ are column vectors.
So the algebra $\mathbf{A}$ (binary operation, bilinear map, tensor) is presented by the matrix $A_e\in Mat(m\times m^2;F)$-the matrix of structure constants (MSC) of $\mathbf{A}$ with respect to the basis $e$.

If $e'=(e'_1,e'_2,...,e'_m)$ is also a basis for $\mathbf{A}$, $g\in GL(m,F)$, $e'g=e$ then it is well known that
\[A_{e'}=gA_e(g^{-1})^{\otimes 2}\] is valid.
Further a basis $e$ is fixed and therefore instead of $A_e$ we use $A$, we do not make difference between $\mathbf{A}$ and its matrix $A$. Let $X=(X^i_{jk})_{i,j,k=1,2,...,m}$ stand for a variable matrix and $Tr_1(X)$, $Tr_2(X)$ stand for the row vectors \[(\sum_{i=1}^mX^i_{i1},\sum_{i=1}^mX^i_{i2},...,\sum_{i=1}^mX^i_{im}),\ \ (\sum_{i=1}^mX^i_{1i},\sum_{i=1}^mX^i_{2i},...,\sum_{i=1}^mX^i_{mi}),\] respectively.

We use $\tau$  for the representation of $GL(m,F)$ on the $n=m^3$ dimensional vector space $V=Mat(m\times m^2;F)$ defined by \[ \tau: (g,A)\mapsto B=gA(g^{-1}\otimes g^{-1}).\]

For simplicity instead of "$\tau$-equivalent", "$\tau$-invariant" we use "equivalent" and "invariant".

We represent each MSC $A$ as a row vector with entries from $Mat(m,F)$ by parting it consequently into elements of $Mat(m,F)$:  \[A=(A_1,A_2,...,A_m),\ A_1,A_2,...,A_m\in Mat(m,F).\]

If $C$ is a block matrix with blocks from $Mat(m,F)$ we use notation $C^{\overline{*}}$, where $*$ is the tensor product or transpose operation, to mean that the operation $*$ with $C$ is done "over $Mat(m,F)$" (not over $F$), for example for the above presented matrix $A$ \[A^{\overline{t}}=\left(\begin{array}{c}
A_1\\ A_2\\ \vdots\\ A_m
\end{array}\right)-\ \ \mbox{column vector over}\ Mat(m,F),\ \]
\[A^{\overline{\otimes} 2}=
(A^2_1, A_1A_2,..., A_1A_m, A_2A_1, A^2_2,..., A_2A_m,..., A_mA_1, A_mA_2,..., A^2_m).\]

 One can see that the equality $B=gA(g^{-1}\otimes g^{-1})$ can be presented as \[B=(B_1,B_2,...,B_m)=gA(g^{-1})^{\otimes 2}=(gA_1g^{-1}, gA_2g^{-1},..., gA_mg^{-1})(g^{-1}\otimes I),\] where $I$ stands for $m\times m$ size identity matrix. Moreover for any matrices $C$ and $D$ the equality \[ (C\otimes I)\overline{\otimes}(D\otimes I)=(C\otimes D)\otimes I\] holds true. Therefore the following equalities hold true.
 \[(B_1, B_2,..., B_m)^{\overline{\otimes} k}=(gA_1g^{-1}, gA_2g^{-1},..., gA_mg^{-1})^{\overline{\otimes} k}((g^{-1})^{\otimes k}\otimes I),\]
 \[\left(\begin{array}{c}
 B_1\\ B_2\\ \vdots\\ B_m
 \end{array}\right)^{\overline{\otimes} k}(B_1, B_2,..., B_m)^{\overline{\otimes} k}=\left(\begin{array}{cccc}
 B_1^2& B_1B_2& \cdots & B_1B_m\\
 B_2B_1& B_2^2& \cdots & B_2B_m\\
 \vdots& \vdots&\cdots&\vdots\\
 B_mB_1&B_mB_2& \cdots & B_m^2
 \end{array}\right)^{\overline{\otimes} k},\]  
 \[\left(\begin{array}{cccc}
 B_1^2& B_1B_2& \cdots & B_1B_m\\
 B_2B_1& B_2^2& \cdots & B_2B_m\\
 \vdots& \vdots&\cdots&\vdots\\
 B_mB_1&B_mB_2& \cdots & B_m^2
 \end{array}\right)^{\overline{\otimes} k}=\] \[ (((g^t)^{-1})^{\otimes k}\otimes I)\left(\begin{array}{cccc}
 gA_1^2g^{-1}& gA_1A_2g^{-1}& \cdots &gA_1A_mg^{-1}\\
 gA_2A_1g^{-1}& gA_2^2g^{-1}& \cdots & gA_2A_mg^{-1}\\
 \vdots& \vdots&\cdots&\vdots\\
 gA_mA_1g^{-1}&gA_mA_2g^{-1}& \cdots & gA_m^2g^{-1}
 \end{array}\right)^{\overline{\otimes} k}((g^{-1})^{\otimes k}\otimes I).\] 
 
 Component-wise application of trace to this equality, which is denoted by $\tilde{Tr}$ results in 
 \[\tilde{Tr}(\left(\begin{array}{cccc}
 B_1^2& B_1B_2& \cdots & B_1B_m\\
 B_2B_1& B_2^2& \cdots & B_2B_m\\
 \vdots& \vdots&\cdots&\vdots\\
 B_mB_1&B_mB_2& \cdots & B_m^2
 \end{array}\right)^{\overline{\otimes} k})=\] \[ ((g^t)^{-1})^{\otimes k}\tilde{Tr}(\left(\begin{array}{cccc}
 gA_1^2g^{-1}& gA_1A_2g^{-1}& \cdots &gA_1A_mg^{-1}\\
 gA_2A_1g^{-1}& gA_2^2g^{-1}& \cdots & gA_2A_mg^{-1}\\
 \vdots& \vdots&\cdots&\vdots\\
 gA_mA_1g^{-1}&gA_mA_2g^{-1}& \cdots & gA_m^2g^{-1}
 \end{array}\right)^{\overline{\otimes} k})(g^{-1})^{\otimes k}=\]
 \[ ((g^{-1})^{\otimes k})^t\tilde{Tr}(\left(\begin{array}{cccc}
 A_1^2& A_1A_2& \cdots &A_1A_m\\
 A_2A_1& A_2^2& \cdots & A_2A_m\\
 \vdots& \vdots&\cdots&\vdots\\
 A_mA_1&A_mA_2& \cdots & A_m^2
 \end{array}\right)^{\overline{\otimes} k})(g^{-1})^{\otimes k},\] as far as for any matrices $C$, $D$ and $E$, where $D$ is a block matrix with blocks from $Mat(m,F)$ and $(C\otimes I)D(E\otimes I)$ has meaning, the equality \[\tilde{Tr}((C\otimes I)D(E\otimes I))=C\tilde{Tr}(D)E\] is valid. One can represent the above obtained matrix equality in the following compact form
 \[\tilde{Tr}((B^{\overline{t}}B)^{\overline{\otimes} k})=((g^{-1})^{\otimes k})^t\tilde{Tr}((A^{\overline{t}}A)^{\overline{\otimes} k})(g^{-1})^{\otimes k}.\] Note that $\tilde{Tr}((A^{\overline{t}}A)^{\overline{\otimes} k})$ is a symmetric matrix. The obtained equality allows formulation of the following theorem.
 
 {\bf Theorem 1.} \textit{Invariants of the quadratic forms given by the matrix $\tilde{Tr}((X^{\overline{t}}X)^{\overline{\otimes} k})$ are invariants of the $m$-dimensional algebras.}
 
 This result can be used for a rough classification of finite dimensional algebras: Two $m$-dimensional algebras $A$, $B$ are rough equivalent if the quadratic forms given by matrices \[\tilde{Tr}(A^{\overline{t}}A)=\left(\begin{array}{cccc}
 Tr(A_1^2)& Tr(A_1A_2)& \cdots & Tr(A_1A_m)\\
 Tr(A_2A_1)& Tr(A_2^2)& \cdots & Tr(A_2A_m)\\
 \vdots& \vdots&\cdots&\vdots\\
 Tr(A_mA_1)&Tr(A_mA_2)& \cdots & Tr(A_m^2)
 \end{array}\right),\] 
 \[\tilde{Tr}(B^{\overline{t}}B)=\left(\begin{array}{cccc}
 Tr(B_1^2)& Tr(B_1B_2)& \cdots & Tr(B_1B_m)\\
 Tr(B_2B_1)& Tr(B_2^2)& \cdots & Tr(B_2B_m)\\
 \vdots& \vdots&\cdots&\vdots\\
 Tr(B_mB_1)&Tr(B_mB_2)& \cdots & Tr(B_m^2)
 \end{array}\right)\] are equivalent.
 
 It is clear that entries of $\tilde{Tr}(X^{\overline{t}}X)$ are polynomials in components of $X$ and there exists nonsingular matrix $Q(X^{\overline{t}}X)$ with rational entries in $X$ such that the matrix \[\tilde{Tr}(\overline{X}^{\overline{t}}\overline{X})=(Q(X^{\overline{t}}X)^{-1})^t\tilde{Tr}(X^{\overline{t}}X)Q(X^{\overline{t}}X)^{-1}=D(X)\] is a diagonal matrix and $Q(g)=I$ whenever $g$ is a nonsingular diagonal matrix, where $\overline{X}=\tau(Q(X^{\overline{t}}X),X)$. 

In algebraically closed field $F$ case it means that one can define a nonempty invariant open subset $V_0\subset V$ such that $\tilde{Tr}(\overline{A}^{\overline{t}}\overline{A})=D(A)$ and $D(A)$ is nonsingular whenever $A\in V_0$.

{\bf Theorem 2.} \textit{ Two algebras $A, B\in V_0$ are equivalent(isomorphic) if and only if 
\[\overline{B}=\tau(g_0,\overline{A})\ \mbox{for some}\ g_0\in GL(m,F)\ \mbox{for which}\ g_0^tD(B)g_0=D(A).\]}
{\bf Proof.} If $B=\tau(g,A)$ then $\overline{B}=\tau(Q(B^{\overline{t}}B),B)=
\tau(Q(B^{\overline{t}}B),\tau(g,A))=$ \[\tau(Q(B^{\overline{t}}B)g,A)=\tau(Q(B^{\overline{t}}B)gQ(A^{\overline{t}}A)^{-1},\tau(Q(A^{\overline{t}}A),A)=\tau(Q(B^{\overline{t}}B)gQ(A^{\overline{t}}A)^{-1},\overline{A}),\] and for $g_0=Q(B^{\overline{t}}B)gQ(A^{\overline{t}}A)^{-1}$ one has 
\[g_0^tD(B)g_0=  (Q(B^{\overline{t}}B)gQ(A^{\overline{t}}A)^{-1})^tD(B)Q(B^{\overline{t}}B)gQ(A^{\overline{t}}A)^{-1}=\]
\[(Q(A^{\overline{t}}A)^{-1})^t(g^t(Q(B^{\overline{t}}B)^tD(B)Q(B^{\overline{t}}B))g)Q(A^{\overline{t}}A)^{-1}=\] \[(Q(A^{\overline{t}}A)^{-1})^t(g^t(\tilde{Tr}(B^{\overline{t}}B))g)Q(A^{\overline{t}}A)^{-1}= (Q(A^{\overline{t}}A)^{-1})^t\tilde{Tr}(A^{\overline{t}}A)Q(A^{\overline{t}}A)^{-1}=D(A).\]

Visa versa if $\overline{B}=\tau(g_0,\overline{A})$ for some $g_0$ for which $g_0^tD(B)g_0=D(A)$ then for\\ $g=Q(B^{\overline{t}}B)^{-1}g_0Q(A^{\overline{t}}A)$ one has 
\[\tau(g,A)=\tau(Q(B^{\overline{t}}B)^{-1}g_0Q(A^{\overline{t}}A),A)= \tau(Q(B^{\overline{t}}B)^{-1}g_0,\tau(Q(A^{\overline{t}}A),A)=\] \[\tau(Q(B^{\overline{t}}B)^{-1}g_0,\overline{A})=\tau(Q(B^{\overline{t}}B)^{-1},\tau(g_0,\overline{A}))=\tau(Q(B^{\overline{t}}B)^{-1},\overline{B})=B.\]

Assume that there exists  matrix $P(X)$, with rational entries with respect to components of $X$, such that $P(\overline{A})$ is nonsingular for any $A\in V_0$ and the equality \begin{equation} P(\tau(g,\overline{A}))=P(\overline{A})g^{-1}\ \mbox{holds true whenever} \ g^tD(\tau(g,A))g=D(A).
\end{equation} 

{\bf Theorem 3.} \textit{ For $A$,$B\in V_0$ there exists $g_0\in GL(m,F)$ such that  $g_0^tD(B)g_0=D(A)$ and $\overline{B}=\tau(g_0,\overline{A})$ if and only if \[\tau(P(\overline{B}),\overline{B})= \tau(P(\overline{A}),\overline{A}),\ (P(\overline{B})^{-1})^tD(B) P(\overline{B})^{-1}=(P(\overline{A})^{-1})^tD(A) P(\overline{A})^{-1}.\]}
{\bf Proof.} If $\overline{B}=\tau(g_0,\overline{A})$ and  $g_0^tD(B)g_0=D(A)$ then 
\[\tau(P(\overline{B}),\overline{B})=\tau(P(\tau(g_0,\overline{A})),\tau(g_0,\overline{A}))=\tau(P(\overline{A})g_0^{-1},\tau(g_0,\overline{A}))=\tau(P(\overline{A}),\overline{A})\] and $(P(\overline{B})^{-1})^tD(B) P(\overline{B})^{-1}=((P(\overline{A})g_0^{-1})^{-1})^tD(B) (P(\overline{A})g_0^{-1})^{-1}=$
\[((P(\overline{A})^{-1})^tg_0^tD(B)g_0 P(\overline{A})^{-1}=(P(\overline{A})^{-1})^tD(A) P(\overline{A})^{-1}.\]

Visa versa, if equalities \[\tau(P(\overline{B}),\overline{B})= \tau(P(\overline{A}),\overline{A}),\ (P(\overline{B})^{-1})^tD(B) P(\overline{B})^{-1}=(P(\overline{A})^{-1})^tD(A) P(\overline{A})^{-1}\] are valid then for $g_0=P(\overline{B})^{-1}P(\overline{A})$ one has
$g_0^tD(B)g_0=D(A)$ and \[\tau(g_0,\overline{A})=\tau(P(\overline{B})^{-1} P(\overline{A}),\overline{A})=\tau(P(\overline{B})^{-1}, \tau(P(\overline{A}),\overline{A}))=\tau(P(\overline{B})^{-1}, \tau(P(\overline{B}),\overline{B}))=\overline{B}.\]

So Theorems 2 and 3 imply that the system of entries of matrices \[\tau(P(\overline{X}),\overline{X}),\ \ (P(\overline{X})^{-1})^t\tilde{Tr}(\overline{X}^{\overline{t}}\overline{X}) P(\overline{X})^{-1}\] is a separating system of rational invariants for algebras from $V_0$.

The above presented results show importance of construction of matrix $P(X)$ with properties (1). Further we  discuss a construction of such matrix by the use of rows $r(\overline{A})$ for which the equality  
\[r(\tau(g,\overline{A}))=r(\overline{A})g^{-1}\] is valid, whenever $g^tD(\tau(g,A))g=D(A)$. To construct such rows one can use the following approach. 

Assume that the equalities \[
	\overline{B}=g\overline{A}(g^{-1})^{\otimes 2},\ 
	\tilde{C}=gCg^t\] are true, where $C^t=C$ and $C$ is a nonsingular matrix.
In this case
\[\tilde{C}^{\otimes 2}=g^{\otimes 2}C^{\otimes 2}(g^{\otimes 2})^t,\  
 \overline{B}\tilde{C}^{\otimes 2}=g\overline{A}C^{\otimes 2}(g^{\otimes 2})^t,\ \tilde{C}^{\otimes 2}\overline{B}^t=g^{\otimes 2}C^{\otimes 2}\overline{A}^tg^t\] and \[\overline{B}\tilde{C}^{\otimes 2}\overline{B}^t=g\overline{A}C^{\otimes 2}\overline{A}^tg^t. \]
 
 On induction it is easy to see that for any natural $k$ the equality
\[\overline{B}^{\otimes 2^0}\overline{B}^{\otimes 2^1}...\overline{B}^{\otimes 2^{k-1}}\tilde{C}^{\otimes 2^k}(\overline{B}^{\otimes 2^0}\overline{B}^{\otimes 2^1}...\overline{B}^{\otimes 2^{k-1}})^t=g\overline{A}^{\otimes 2^0}\overline{A}^{\otimes 2^1}...\overline{A}^{\otimes 2^{k-1}}C^{\otimes 2^k}(\overline{A}^{\otimes 2^0}\overline{A}^{\otimes 2^1}...\overline{A}^{\otimes 2^{k-1}})^tg^t \] holds true. Therefore due to the equalities 
\[\overline{B}^{\otimes 2^0}\overline{B}^{\otimes 2^1}...\overline{B}^{\otimes 2^{k-1}}\tilde{C}^{\otimes 2^k}(\overline{B}^{\otimes 2^0}\overline{B}^{\otimes 2^1}...\overline{B}^{\otimes 2^{k-1}})^t\tilde{C}^{-1}=\] \[g\overline{A}^{\otimes 2^0}\overline{A}^{\otimes 2^1}...\overline{A}^{\otimes 2^{k-1}}C^{\otimes 2^k}(\overline{A}^{\otimes 2^0}\overline{A}^{\otimes 2^1}...\overline{A}^{\otimes 2^{k-1}})^t C^{-1}g^{-1}, \]
\[\overline{B}^{\otimes 2^0}\overline{B}^{\otimes 2^1}...\overline{B}^{\otimes 2^{k-1}}\tilde{C}^{\otimes 2^k}(\overline{B}^{\otimes 2^0}\overline{B}^{\otimes 2^1}...\overline{B}^{\otimes 2^{k-1}})^t\tilde{C}^{-1}\overline{B}=\] \[g\overline{A}^{\otimes 2^0}\overline{A}^{\otimes 2^1}...\overline{A}^{\otimes 2^{k-1}}C^{\otimes 2^k}(\overline{A}^{\otimes 2^0}\overline{A}^{\otimes 2^1}...\overline{A}^{\otimes 2^{k-1}})^t C^{-1}\overline{A}(g^{-1})^{\otimes 2} \] one has
\[Tr_i(\overline{B}^{\otimes 2^0}\overline{B}^{\otimes 2^1}...\overline{B}^{\otimes 2^{k-1}}\tilde{C}^{\otimes 2^k}(\overline{B}^{\otimes 2^0}\overline{B}^{\otimes 2^1}...\overline{B}^{\otimes 2^{k-1}})^t\tilde{C}^{-1}\overline{B})=\] \[Tr_i(\overline{A}^{\otimes 2^0}\overline{A}^{\otimes 2^1}...\overline{A}^{\otimes 2^{k-1}}C^{\otimes 2^k}(\overline{A}^{\otimes 2^0}\overline{A}^{\otimes 2^1}...\overline{A}^{\otimes 2^{k-1}})^t C^{-1}\overline{A})g^{-1},\ i=1,2. \]

The last equality shows that in our algebra case on can try to construct the needed matrix $P(X)$ by the use of rows
\[Tr_i(X^{\otimes 2^0}X^{\otimes 2}...X^{\otimes 2(k-1)}((\tilde{Tr}(X^{\overline{t}}X))^{-1})^{\otimes 2^k}(X^{\otimes 2^0}X^{\otimes 2}...X^{\otimes 2(k-1)})^t \tilde{Tr}(X^{\overline{t}}X)X),\]
where $i=1,2,\ k=0,1,2,...$.
 
What is left unjustified here is that one should justify existence, in general, of a linear independent system consisting of $m$ such rows.

{\bf Remark.} \textit{After a rough classification one can classify further each case of the rough classification with respect to the corresponding stabilizer.} 
\begin{theacknowledgments}The author acknowledges MOHE for a support by grant FRGS14-153-0394.\end{theacknowledgments}
\vskip 0.4 true cm
\bibliographystyle{aipproc}

\end{document}